\newtheorem{defn}{Definition}[section]
\newtheorem{lem}[defn]{Lemma}
\newtheorem{thm}[defn]{Theorem}
\newtheorem{prop}[defn]{Proposition}
\newtheorem{cor}[defn]{Corollary}
\newtheorem{rem}[defn]{Remark}

\newtheorem{assu}[defn]{Assumption}

\newtheorem{nota}[defn]{Notation}
 
\newcommand{\C}{{\bf C}}
\newcommand{\PP}{{\bf P}}
\newcommand{\Q}{{\bf Q}}
\newcommand{\Z}{{\bf Z}}
\newcommand{\cover}{f:X\to Y}
\newcommand{\epsi}{\epsilon}
\newcommand{\Ho}{{\rm H}}
\newcommand{\HH}[1]{{\rm H}^{#1}}
\newcommand{\Hom}{{\rm Hom}}
\newcommand{\ichi}{_{\chi}}
\newcommand{\ifi}{_{\phi}}
\newcommand{\ichifi}{_{\chi,\phi}}
\newcommand{\ichichi}{_{\chi,\chi\inv}}
\newcommand{\ichinv}{_{\chi\inv}}
\newcommand{\ifinv}{_{\phi\inv}}

\newcommand{\ichifinv}{_{\chi\phi,(\chi\phi)\inv}}

\newcommand{\achi}{^{(\chi)}}

\newcommand{\inv}{^{-1}}
\newcommand{\F}{{\cal F}}
\newcommand{\I}{{\cal I}}
\newcommand{\OX}{{\cal O}_X}
\newcommand{\OY}{{\cal O}_Y}

\newcommand{\proof}{{\bf Proof:\,}}
\newcommand{\Pic}{{\rm Pic}}
\newcommand{\qed}{\ $\Box$\par\smallskip}
\newcommand{\tensor}{\!\otimes\!}

\title{On the period map for abelian covers of projective varieties}
\author{Rita Pardini - Universit\`a di Udine}
\date{}
\documentstyle{article}
\begin{document}
\maketitle
\setcounter{section}{-1}
\def\theequation{\thesection.\arabic{equation}}

\section{Introduction}
\setcounter{defn}{0}
\setcounter{equation}{0}

This paper is devoted to the study of the period map for abelian covers
of smooth projective varieties of dimension $n\ge 2$. Our viewpoint is
very close to that of Green in \cite{suffampio}, namely we look for 
 results that hold for abelian covers of an arbitrary variety whenever
certain ampleness assumptions on the building data defining the cover
are satisfied. We focus on two questions: the infinitesimal and the
variational Torelli problems.

Infinitesimal Torelli for the periods of $k$-forms holds for a smooth
projective variety
$X$ if
the map  $\Ho^1(X,T_X)\to
\oplus_p\Hom\left (\Ho^p(X,\Omega_X^{k-p}),\Ho^{p+1}(X,
\Omega_X^{k-p-1})\right)$, expressing the differential of the period map
for $k$-forms, is injective. This is expected to be true as soon as the
canonical bundle of
$X$ is ``sufficiently ample''. There are many results in this direction,
concerning special classes of varieties, as hypersurfaces (see, for
instance,
\cite{suffampio}), complete intersections (\cite{konnocompl}) and simple
cyclic covers (\cite{konnociclico},
\cite{lwp},\cite{peters}). Here we continue the work on abelian covers
of  
\cite{torelli}, and prove (see
\ref{mainthm1}):
\begin{thm}\label{introinf} 
Let $G$ be an abelian group and let
$\cover$ be a $G$-cover, with $X$, $Y$
smooth projective varieties of dimension $n\ge2$. If properties
$(A)$ and
$(B)$ of
\ref{ipotesi} are satisfied, then infinitesimal Torelli for the periods
of $n$-forms holds for $X$.
\end{thm}
Properties $(A)$ and $(B)$ amount to the vanishing of certain cohomology
groups and are certainly satisfied if the building data of the cover are
sufficiently ample. If $Y$ is a special variety (e.g.,
$Y=\PP^n$), then thm. \ref{introinf} yields an almost sharp
statement  (see theorem
\ref{mainthm1p}). In general, a result of 
Ein-Lazarsfeld (\cite{el}) and Griffiths vanishing theorem enable us (see
prop.
\ref{effective}) to give explicit conditions under which $(A)$ and $(B)$
are satisfied, and thus to deduce an effective
statement  from thm. \ref{introinf} (see thm. \ref{mainthm1e}).

In order to extend to the case of
arbitrary varieties the infinitesimal Torelli theorem obtained in
\cite{torelli} for a special class of surfaces, we introduce a
generalized notion of prolongation bundle and  give a Jacobi ring
construction analogous to those of \cite{suffampio} and
\cite{konnovar}. This is also a starting point for  attacking the
variational Torelli problem, which asks whether, given a flat family
${\cal X}\to B$ of smooth polarized varieties, the map associating to a
point
$b\in B$ the infinitesimal variation of Hodge structure of the fibre
$X_b$ is generically injective, up to isomorphism of polarized
varieties.   A positive answer to this problem has been given for
families of projective hypersurfaces (\cite{cggh},\cite{donagi}), for
hypersurfaces of high degree of arbitrary varieties (\cite{suffampio}),
for some complete intersections (\cite{konnovar}) and for simple cyclic
covers of high degree (\cite{ciclico}). The most effective tool in
handling these problems is the symmetrizer, introduced by Donagi, but
unfortunately, an  analogous construction does not seem feasible in the
case of abelian covers. However, exploiting the variational Torelli
result of
\cite{suffampio}, we are able to obtain, under analogous assumptions, a
similar result for a large class of  abelian covers. 
More precisely, we prove (see thm. \ref{mainthm2}):
\begin{thm}\label{intro2}
(Notation as in section \ref{covers}.)

\noindent Let $Y$ be a smooth projective variety of dimension $n\ge 2$,
with very ample canonical class. Let $G$
be a finite abelian group and let 
$\cover$ be a smooth $G$-cover with sufficiently ample building data
$L\ichi$, $D_i$, $\chi\in G^*$, $i=1,\ldots r$. Assume that for every
$i=1,\ldots r$ the identity is the only automorphism of
$Y$ that preserves the linear equivalence class of $D_i$; moreover,
assume that for $i=1,\ldots r$ there exist a
$\chi\in G^*$ (possibly depending on $i$) such that $\chi(g_i)\ne 1$ and
$L\ichi(-D_i)$ is ample. Let
${\cal X}\to
\tilde{W}$ be the family of the smooth
$G$-covers of
$Y$ obtained by letting the $D_i$'s vary in their linear equivalence
classes: then there is a dense open set $V\subset\tilde{W}$ such that
the fibre $X_s$ of ${\cal X}$ over $s\in V$ is determined by its IVHS
for $n$-forms plus the natural $G$-action on it.
\end{thm}
The paper is organized as follows:  section $1$ is a brief review of
abelian covers, sections $2$ and $3$ contain the technical details about
prolongation bundles and the Jacobi ring construction, section $4$
contains the statements and proofs of the results on infinitesimal
Torelli, and section
$5$ contains some technical lemmas that allow us to prove in section $6$
the variational Torelli theorem \ref{intro2}.

\noindent {\em Acknowledgements:} I wish to thank Mark Green for
communicating me the proof of lemma \ref{green}.

 \paragraph{}
{\bf Notation and conventions:} All varieties are smooth projective
varieties of dimension $n\ge 2$ over the field $\C$ of complex numbers. We
do not distinguish between vector bundles and locally free sheaves; as a
rule, we use the additive notation for divisors and the multiplicative
notation for line bundles.  For a divisor
$D$, $c_1(D)$ denotes the first Chern class of $D$ and $|D|$ the complete
linear system of $D$. If $L$ is a line bundle, we also denote by $|L|$ the
complete linear system of $L$, and  we write $L^k$ for $L^{\tensor^k}$
and $L\inv$ for the dual line bundle. As usual,
$T_Y$ denotes the tangent sheaf of
$Y$,
$\Omega^k_Y$ denotes the sheaf of regular $k$-forms on $Y$,
$\omega_Y=\Omega^n_Y$ denotes the canonical bundle and $\Pic(Y)$ the
Picard group. 
 If ${\cal F}$ is a locally free sheaf, we denote by ${\cal F}^*$ the dual
sheaf,  by $S^k\F$ the $k$-th symmetric power of $\F$ and by $\det {\cal
F}$ the determinant bundle. Consistently, the dual of a vector space $U$
is denoted by
$U^*$; the group of linear automorphisms of $U$ is denoted by $GL(U)$.

\noindent $[x]$ denotes the integral part of the natural number $x$.

\section{Abelian covers and projection formulas}\label{covers}
\setcounter{defn}{0}
\setcounter{equation}{0}

In this section we recall some facts about abelian covers that will be
needed later. For more details and proofs, see \cite{abeliani}.

Let $G$ be a finite abelian group of order $m$ and let $G^*=\Hom(G,\C^*)$
be the group of characters of $G$.  A {\em $G$-cover} of a smooth
$n$-dimensional variety
$Y$ is a Galois cover $\cover$ with Galois group $G$, with $X$ normal. Let
${\cal F}$ be a
$G$-linearized locally free 
sheaf of $\OX$-modules: under the action of
$G$, the sheaf
$f_*{\cal F}$ splits as the direct sum of the eigensheaves corresponding
to the characters of $G$. We denote by $(f_*{\cal F})^{(\chi)}$ the
eigensheaf corresponding to a character $\chi\in G^*\setminus\{1\}$ and by
$(f_*{\cal F})^{inv}$ the invariant subsheaf. In particular, when
${\cal F}=\OX$ , we have $(f_*\OX)^{inv}=\OY$ and
$(f_*\OX)^{(\chi)}=L\ichi\inv$, with $L\ichi$ a line bundle.  Let
$D_1,\ldots D_r$ be the irreducible components of the branch locus
$D$ of $f$. For each index $i$, 
 the subgroup of $G$ consisting of the elements that fix the inverse
image of $D_i$ pointwise is a cyclic group $H_i$, the so-called {\em
inertia subgroup} of $D_i$.  The order $m_i$ of $H_i$ is equal to the
order of ramification of
$f$ over $D_i$ and the representation of $H_i$ obtained by taking
differentials and restricting to the normal space to
$D_i$ is a faithful character $\chi_i$. 
 The choice of a primitive $m$-th root $\zeta$ of $1$ defines a map from
$\{1,\ldots r\}$ to $G$: the image
$g_i$ of $i$ is the generator of $H_i$ that is mapped to
$\zeta^{m/m_i}$ by $\chi_i$. The line bundles
$L\ichi$, $\chi\in G^*\setminus\{1\}$, and the divisors $D_i$, each
``labelled'' with an element $g_i$ of
$G$ as explained above, are the {\em building data} of the cover, and
determine $\cover$ up to isomorphism commuting with the covering maps.
The building data satisfy the so-called {\em fundamental relations}. In
order to write these down, we have to set some notation. For $i=1,\ldots
r$ and $\chi\in G^*$, we denote by
$a^i\ichi$ the smallest positive integer such that
$\chi(g_i)=\zeta^{ma^i\ichi/m_i}$; for each pair of characters
$\chi$,$\phi$ we set $\epsi\ichifi^i=[(a^i\ichi+a^i\ifi)/m_i]$ (notice
that $\epsi\ichifi^i=0$ or $1$) and
$D\ichifi=\sum_{i=1}^r\epsi\ichifi^iD_i$. In particular, $D\ichichi$ is
the sum of the components $D_i$ of $D$ such that $\chi(g_i)\ne 1$. Then,
the fundamental relations of the cover are the following:
\begin{equation} L\ichi+L\ifi\equiv L_{\chi\phi}+ D\ichifi, \qquad
\forall \chi,\phi\in G^*
\label{fundrel}
\end{equation} When $\phi=\chi\inv$, the fundamental
relations read:
\begin{equation}\label{fundrelbis} L\ichi+L\ichinv\equiv D\ichichi.
\end{equation}
The cover $\cover$ can be reconstructed from the building data as
follows: if one chooses sections $s_i$ of $\OY(D_i)$ vanishing on
$D_i$ for
$i=1,\ldots r$, then $X$ is defined inside the vector bundle
$V=\oplus_{\chi\ne 1}L\ichi$ by the equations:
\begin{equation}\label{equazioni}
z\ichi
z\ifi=\left(\Pi_{i}s_i^{\epsilon^i\ichifi}\right)z_{\chi\phi},
\quad\forall \chi,\phi\in G^*\setminus\{1\}
\end{equation}
where
$z\ichi$ denotes the tautological section of the pull-back of $L\ichi$
to $V$. Conversely, for every choice of the sections $s_i$,  equations
(\ref{equazioni}) define a scheme $X$, flat over $Y$, which is
smooth iff the zero divisors of the
$s_i$'s are smooth, their union has only normal crossings singularities
and, whenever $s_{i_1},\ldots s_{i_t}$ all vanish at a point $y$ of $Y$,
the group
$H_{i_1}\times\cdots\times H_{i_t}$ injects into $G$.  So, by letting
$s_i$ vary in $\HH{0}(Y,\OY(D_i))$, one obtains a flat family ${\cal X}$
of smooth
$G$-covers of
$X$, parametrized by an open set $W\subset\oplus_i\HH{0}(Y,\OY(D_i))$.

 Throughout all the paper we will make the following 
\begin{assu}\label{ampleness}
The $G$-cover $\cover$ is smooth of dimension $n\ge 2$; the building data
$L\ichi$,
$D_i$ and the adjoint bundles $\omega_Y\tensor L\ichi$,
$\omega_Y(D_i)$ are ample for every $\chi\in G^*\setminus\{1\}$ and for
every
$i=1,\ldots r$.
\end{assu}
 Assumption \ref{ampleness} implies that the cover is {\em totally
ramified}, namely that
$g_1,\ldots g_r$ generate $G$. Actually, this  is equivalent to the fact
that the divisor $D\ichichi$ is nonempty if the character $\chi$ is
nontrivial, and also to the fact that none of the line bundles $L\ichi$,
$\chi\in G^*\setminus\{1\}$, is a torsion point in $\Pic(Y)$.
Since $X$ is smooth, assumption
\ref{ampleness} implies in particular that for each subset $\{i_1,\ldots
i_t\}\subset\{1,\ldots r\}$, with
$t\leq n$, the cyclic subgroups generated by $g_{i_1},\ldots g_{i_t}$ give
a direct sum inside $G$.

\paragraph{} In principle, all the geometry of $X$ can be recovered from
the geometry of $Y$ and from the building data of
$\cover$. The following proposition is an instance of this philosophy.
\begin{prop}\label{projform} Let $\cover$ be a $G$-cover, with $X$, $Y$
smooth of dimension $n$. For
$\chi\in G^*$, denote by
$\Delta\ichi$ the sum of the components $D_i$ of $D$ such that
$a\ichi^i\neq m_i-1$. Then, for $1\le k\le n$ there are natural
isomorphisms:
$$(f_*\Omega_X^k)\achi=\Omega^k_Y(\log D\ichichi)\tensor L\ichi\inv$$
$$(f_*T_X)\achi=T_Y(-\log\Delta\ichi)\tensor L\ichi\inv$$ and,in
particular:
$$(f_*\Omega^k_X)^{inv}=\Omega^k_Y, \quad (f_*T_X)^{inv}=T_Y(-\log D),
\quad (f_*\omega_X)\achi=\omega_Y\tensor L\ichinv.$$

\end{prop}
\proof This is a slight generalization of Proposition $4.1$ of
\cite{abeliani}, and it can be proven along the same lines. The
identification
$(f_*\omega_X)\achi=\omega_Y\tensor L\ichinv$ follows from the general
formula and relations (\ref{fundrelbis}).\qed

We recall the following generalized form of Kodaira vanishing (see
\cite{EV}, page 56):
\begin{thm}\label{kodaira} Let $Y$ be a smooth projective $n$-dimensional
variety and let $L$ be an ample line bundle. Then:
$$\Ho^i(Y,\Omega^k_Y\tensor L\inv)=0, \quad i+k<n.$$ Moreover, if $A+B$
is a reduced effective normal crossing divisor, then:
$$\Ho^i(Y,\Omega^k_Y(\log(A+B)\tensor L\inv(-B))=0,  \quad i+k<n.$$
\end{thm}

From prop.\ref{projform}, theorem \ref{kodaira} and assumption
\ref{ampleness} it follows that the non-invariant part of the
cohomology of $X$ is concentrated in dimension
$n$. Thus we will be   concerned only with the period map for the
periods of $n$-forms.

\section{Logarithmic forms and sheaf resolutions}\label{logarithmic}
\setcounter{defn}{0}
\setcounter{equation}{0}

In this section we recall the definition and some properties of
logarithmic forms and introduce a generalized notion of prolongation
bundle. We would like to mention that Konno, when studying in
\cite{konnovar} the global Torelli problem for complete intersections,
has also introduced a generalization of the definition of prolongation
bundle, which is however different from the one used here. 

 Let
$D$ be a normal crossing divisor with smooth components
$D_1,\ldots D_r$ on the smooth $n$-dimensional variety $Y$. As usual, we
denote by
$\Omega_Y^k(\log D)$ the sheaf of
$k$-forms with at most logarithmic poles along $D_1,\ldots D_r$ and by
$T_Y(-\log D)$ the subsheaf of $T_Y$ consisting of the vector fields
tangent to the components of $D$. Assume that $y\in Y$ lies
precisely on the components
$D_1,\ldots D_t$ of $D$, with $t\le n$.  Let $x^1,\ldots x^t$ be local
equations for $D_1, \ldots D_t$ and choose $x_{t+1},\ldots x_n$ such that
$x_1,\ldots x_n$ are a set of parame\-tres at $y$. Then
$\frac{dx_1}{x_1},\ldots \frac{dx_t}{x_t}$, $dx_{t+1},\ldots dx_n$ are
a set of free gen\-er\-a\-tors for
$\Omega^1_Y(\log D)$  and
$x_1\frac{\partial}{\partial x_1},\ldots x_t\frac{ \partial}{\partial
x_t}$, $\frac{\partial }{\partial x_{t+1}},\ldots\frac
{\partial}{\partial x_n}$ are free generators for $T_Y(-\log D)$
in a neighbourhood of $y$. So the sheaves of logarithmic forms are
locally free and one has the following canonical identifications: 
$\Omega_Y^k(\log D)=\wedge^k\Omega_Y^1(\log D)$ and
$T_Y(-\log D)=\Omega_Y^1(\log D)^*$, duality being given by contraction of
tensors.  Moreover, we recall that, if ${\cal F}$ is a locally
free sheaf of rank $m$ on $Y$, then the alternation map $\wedge^i{\cal
F}\tensor\wedge ^{m-i}{\cal F}\to
\det{\cal F}$ is a nondegenerate pairing, which induces a canonical
isomorphism
$\left(\wedge^i{\cal F}\right)^*\to \wedge^{m-i}{\cal F}\tensor(\det {\cal
F})\inv$. So we have:
\begin{eqnarray}\label{dualita'}
 T_Y(-\log D)\cong \Omega^{n-1}_Y(\log D)\tensor (\omega_Y(D))\inv \\
\Omega_Y^k(\log D)^*\cong \Omega_Y^{n-k}(\log D)\tensor (\omega_Y(D))\inv
\nonumber
\end{eqnarray} 

\paragraph{} The {\em (generalized) prolongation bundle} $P$ of
$(D_1,\ldots D_r)$ is defined as the ex\-ten\-sion $0\to\Omega^1_Y\to
P\to\oplus^r\OY\to 0$ as\-so\-cia\-ted to the class
$(c_1(D_1),\ldots c_1(D_r))$ of $\Ho^1(Y,\oplus^r\Omega^1_Y)$. Let
$\{U_{\alpha}\}$ be a finite affine covering of $Y$, let $x^i_{\alpha}$ be
local equations for $D_i$ on $U_{\alpha}$, $i=1\ldots r$, and let
$g^i_{\alpha\beta}=x^i_{\alpha}/x^i_{\beta}$. 
 Denote by
$e^1_{\alpha},\ldots e^r_{\alpha}$ the standard basis of
$\oplus^r\OY|_{U_{\alpha}}$.

The elements of
$P|_{U_{\alpha}}$ are represented  by pairs 
$(\sigma_{\alpha}, \sum_i z^i_{\alpha}e_i)$, where
$\sigma_{\alpha}$ is a
$1$-form and the
$z^i_{\alpha}$'s are regular functions, satisfying the following
transition relations  on
$U_{\alpha}\cap U_{\beta}$: 
$$\left(\sigma_{\alpha}, \sum_i
z^i_{\alpha}e^i_{\alpha}\right)=\left(\sigma_{\beta}+\sum_i
z^i_{\beta}\frac{dg^i_{\alpha\beta}}{g^i_{\alpha\beta}},\,
\sum_iz^i_{\beta}e^i_{\beta}\right).$$ There is a natural short exact
sequence:
\begin{equation}\label{olog} 0\to\oplus_i\OY(-D_i)\to P\to\Omega^1_Y(\log
D)\to 0
\end{equation} with dual sequence:
\begin{equation}\label{tlog} 0\to T_Y(-\log D)\to P^*\to\oplus_i
\OY(D_i)\to 0.
\end{equation} In local coordinates the map $\oplus_i
\OY(-D_i)\to P$ is defined by: $x^i_{\alpha}\mapsto
(dx^i_{\alpha},x^i_{\alpha}e^i_{\alpha})$ and the map  $P\to
\Omega^1(\log D)$ is defined by: $(\sigma_{\alpha}, \sum_i
z^i_{\alpha}e^i_{\alpha})\mapsto\sigma_{\alpha}-\sum
z^i_{\alpha}dx^i_{\alpha}/x^i_{\alpha}$.

We close this section by writing down a resolution of
the sheaves of logarithmic forms that will be used in section \ref{ivhs}.
 Given an exact sequence
$0\to A\to B\to C\to 0$ of locally free sheaves, for any $k\ge 1$ one has
the following long exact sequence (see
\cite{cime}, page 39):
\begin{equation}\label{complex}
0\to S^kA\to B\tensor
S^{k-1}A\to\ldots\to
\wedge^{k-1}B\tensor A\to\wedge^k B\to \wedge^kC\to 0.
\end{equation}
Applying this
to (\ref{olog}) and setting $V=\oplus_i\OY (D_i)$ yields:
\begin{equation}\label{reslog}
 0\to S^kV^*\to S^{k-1}V^*\tensor P\to\ldots
\to V^*\tensor\wedge^{k-1}P\to\wedge^kP\to \Omega^k_Y(\log D)\to 0
\end{equation}

\section{The algebraic part of the IVHS}\label{ivhs}
\setcounter{defn}{0}
\setcounter{equation}{0}

 The aim of this section is to give, in the case of abelian covers, a
construction analogous to the Jacobian ring construction for
hypersurfaces of
\cite{suffampio}.

Let ${\cal X}\to B$ be a flat family of smooth projective varieties of
dimension $n$ and let $X$ be the fibre of ${\cal X}$ over the point $0\in
B$; the differential of the period map  for the periods of
$k$-forms for ${\cal X}$ at $0$ is the composition of the Kodaira-Spencer
map with the following universal map, induced by cup-product:
\begin{equation}
\Ho^1(X,T_X)\to \oplus_p\Hom (\Ho^p(X,\Omega_X^{k-p}),\Ho^{p+1}(X,
\Omega_X^{k-p-1})),
\label{inftor}
\end{equation}
This map is called the {\em algebraic part of the infinitesimal
variation of Hodge structure} of
$X$ (IVHS for short). 
Assume that $\cover$ is a smooth $G$-cover; then the
$G$-action on the tangent sheaf and on the sheaves of differential forms
is compatible with cup-product, so the map (\ref{inftor}) splits as the
direct sum of the maps
\begin{equation}\label{inftorab}
\rho^k\ichifi :\Ho^1(X,T_X)^{(\chi)}\to
\oplus_p\Hom(\Ho^p(X,\Omega^{k-p}_X)^{(\phi)},\Ho^{p+1}(X,
\Omega_X^{k-p-1})^{(\chi\phi)})
\end{equation}
 As we have remarked at the end of section
\ref{covers}, if $\cover$ satisfies the assumption \ref{ampleness}, then
the non invariant part of the Hodge structure is concentrated in the
middle dimension $n$, and so we will only describe the IVHS for $k=n$.

We use the
notation of section
\ref{covers} and moreover,
 in order to keep  formulas readable, we set: 
$$T^{inv}=\HH{1}(Y,T_Y(-\log D));\quad
U^{k,inv}=\HH{k}(Y,\Omega_Y^{n-k})$$
$$U^{k,\chi}=\HH{k}(Y,\Omega^{n-k}_Y(\log D\ichichi)\tensor
L\ichi\inv),\qquad k=0,\ldots n,\quad \chi\in G^*\setminus\{1\}.$$
Given a character $\chi\in G^*$,  let  $D_{i_1},\ldots D_{i_s}$ be the
components of $D\ichichi$; let $P^{\chi}$ be the generalized prolongation
bundle of\, $(D_{i_1},\ldots D_{i_s})$ (see section \ref{logarithmic}) and
let
$V^{\chi}=\oplus_j\OY(D_{i_j})$.

Consider the map $(P^{\chi})^*\to V^{\chi}$ defined in sequence
(\ref{tlog}); tensoring this map with $S^{k-1}(V^{\chi})$ and composing
with the symmetrization map $S^{k-1}(V^{\chi})\tensor V^{\chi}\to
S^k(V^{\chi})$, one obtains a map:
\begin{equation}\label{simm} S^{k-1}(V^{\chi})\tensor (P^{\chi})^*\to
S^k(V^{\chi})
\end{equation} Given a line bundle $L$ on $Y$, we define $R^{k,\chi}_L$
to be the cokernel of the map:
$$\HH{0}(Y,S^{k-1}(V^{\chi})\tensor (P^{\chi})^*\tensor
L)\to\HH{0}(Y,S^k(V^{\chi})\tensor L),$$ obtained from (\ref{simm}) by
tensoring with $L$ and passing to global sections. We set
$R^{\chi}_L=\oplus_{k\ge 0}R^{k,\chi}_L$;
for $L=\OY$, $R^{\chi}=R^{\chi}_{\OY}$ is a graded ring and, in
general $R^{\chi}_L$ is a module over $R^{\chi}$, that
we call the {\em Jacobi module} of $L$.  Moreover, if
$L_1$ and
$L_2$ are line bundles on $Y$, then there is an obvious multiplicative
structure:
$$R^{k,\chi}_{L_1}\tensor R^{h,\chi}_{L_2}\to R^{k+h,\chi}_{L_1\tensor
L_2}.$$
 In order to  establish the relationship between the Jacobi
modules and the IVHS of the cover $X$, we need some definitions.

\begin{defn}\label{ipotesi}
For a $G$-cover $\cover$ satisfying assumption \ref{ampleness}, let
$\Gamma_f$ be the semigroup of $\Pic(Y)$ generated by the building data.
We say that:

\noindent $X$ has property $(A)$ iff $\HH{k}(Y,\Omega_Y^j\tensor
L)=0$ and $\HH{k}(Y,\Omega_Y^j\tensor\omega_Y\tensor L)=0$ for
$k>0$, $j\ge 0$, 
$L\in\Gamma_f\setminus\{0\}$;

\noindent $X$ has property $(B)$ iff for $L_1$, $L_2$ in
$\Gamma_f\setminus\{0\}$ the  mul\-ti\-pli\-ca\-tion map
$\HH{0}(Y,\omega_Y\tensor L_1)\otimes\HH{0}(Y,\omega_Y\tensor L_2)\to
\HH{0}(Y,\omega_Y^2\tensor L_1\tensor L_2)$ is surjective.
\end{defn}
\begin{rem}\label{annullamento}
 If $M$ is an ample line bundle such that 
$\HH{k}(Y,\Omega_Y^j\tensor M)=0$ for $k>0$ and $j\ge 0$, then the
cohomology groups $\HH{k}(Y,\wedge^jP\tensor M\inv)$ vanish for $j\ge 0$
and
$k<n$.
\end{rem}
\proof By Serre duality, it is equivalent to show that
$\HH{r}(Y,\wedge^jP^*\tensor M\tensor \omega_Y)=0$ for $r>0$. In turn,
this can be proven by induction on $j$, by looking at the
hypercohomology of the complex obtained by applying (\ref{complex}) to
the sequence $0\to \oplus_i\OY\to P^*\to T_Y\to 0$.
\qed

 We also introduce the following 
\begin{nota}
Let $L$ and $M$ be line bundles on the smooth variety $Y$; if $L\otimes
M\inv$ is ample, then we  write $L>M$ and, if $L\tensor M\inv$ is nef,
then we  write $L\ge M$. We use the same notation for divisors.
\end{nota}
\begin{rem}\label{projective}
Properties $(A)$ and $(B)$ are easily checked
for coverings of certain varieties $Y$; for instance, if\/ $Y=\PP^n$, then
by Bott vanishing theorem it is enough to require that $\omega_Y\tensor
L\ichi>0$ and
$\omega_Y(D_i)>0$ for every $\chi\ne 1$ and for $i=1,\ldots r$.
\end{rem}
The
next proposition yields an effective criterion for $(A)$ and $(B)$ in
case $Y$ is an arbitrary variety.
\begin{prop}\label{effective}
Let $\cover$ be a $G$-cover satisfying assumption \ref{ampleness} and let
$E$ be a very ample divisor on $Y$. Define
$c(n)=\left(\!\!\begin{array}{c} n-1\\
(n-1)/2\end{array}\!\!\right)$ if
$n$ is odd and
$c(n)=
\left(\!\!\begin{array}{c} n-1 \\ n/2\end{array}\!\!\right)$ if
$n$ is even, and set $E_n=\left(\omega_Y(2nE)\right)^{c(n)}$. 

i) if $D_i$, $L\ichi$, $\omega_Y(D_i)$, $\omega_Y\tensor
L\ichi>E_n$ for $\chi\ne 1$ and for $i=1,\ldots r$, then $(A)$ is
satisfied.

ii) if $L\ichi$ and $D_i\ge (n+1)E$  for $\chi\ne 1$ and for $i=1,\ldots
r$, then
$(B)$ is satisfied.
\end{prop}
\proof
The complete linear system $|E|$ embeds $Y$ in a projective space $\PP$;
since $\Omega_{\PP}^{j}(j+1)$ is generated by global sections, the sheaf
$W^j=\Omega_Y^j((j+1)E)$, being a quotient of the former bundle, is also
generated by global sections. By Griffiths vanishing theorem (\cite{SS},
theorem 5.52), if
$N$ is an ample line bundle, then the cohomology group
$\HH{k}(Y,W^j\tensor\omega_Y\tensor
\det(W^j)\tensor N)$ vanishes for $k>0$. We recall from adjunction
theory that  $\omega_Y((n+1)E)$ is base point free and therefore nef;
using this fact, it is easy to check that $E_n\ge\det(W^j)$ for $j\ge 0$.
Statement i) now follows immediately from Griffiths vanishing. 

In order to prove ii), set $V_1=\HH{0}(Y,\omega_Y\tensor L_1)$. The
assumptions imply that
$V_1$ is base-point free, so evaluation of sections gives the following
short exact sequence of locally free sheaves:
$0\to K_1\to\OY\tensor V_1\to \omega_Y\tensor L_1\to 0$. Twisting with
$\omega_Y\tensor L_2$ and passing to cohomology, one sees that the
statement follows if
$\HH{1}(Y,K_1\tensor\omega_Y\tensor L_2)=0$. In turn, this is precisely
case
$k=q=1$ of theorem $2.1$ of \cite{el}.
\qed
\begin{lem}\label{iso}
Let $\cover$ be a $G$-cover satisfying property $(A)$ and let $L$ be a
line bundle on $Y$; if $L$ or $L\tensor\omega_Y\inv$ belong to
$\Gamma_f\setminus\{0\}$, then for every $\chi\in G^*\setminus\{1\}$ there
is a natural isomorphism:
$$\HH{k}(Y,\Omega_Y^{n-k}(\log D\ichichi)\tensor L\inv)\cong
(R^{n-k,\chi}_{\omega_Y\tensor L})^*$$
In particular, there are natural isomorphisms:
 $$T^{inv}\cong(R^{n-1,1}_{\omega_Y^2(D)})^*;\qquad U^{k,\chi}\cong
(R^{n-k,\chi}_{\omega_Y\tensor L\ichi})^*,\quad \chi\ne 1 $$
$$\HH{1}\left(Y,\Omega^{n-1}_Y(\log D\ichichi)\tensor(L\ichi\tensor
L\ifinv\tensor\omega_Y)\inv\right)\cong (R^{n-1,\chi}_{\omega_Y^2\tensor
L\ichi\tensor L\ifinv})^*,\quad\chi\ne 1.$$
\end{lem}
\proof
For $k=n$, the statement is just Serre duality. For $k<n$, we compute
$\HH{k}\left(Y,\Omega_Y^{n-k}(\log D\ichichi)\tensor L\inv\right)$ by
ten\-sor\-ing the res\-o\-lu\-tion (\ref{reslog}) of the sheaf
$\Omega_Y^{n-k}(\log D\ichichi)$  with
$L\inv$ and breaking up the resolution thus
obtained into short exact sequences. Remark \ref{annullamento} and theorem
\ref{kodaira} imply that the cohomology groups
$\HH{n-k+j}(Y,S^j(V\ichi)^*\tensor L\inv\tensor \wedge^{n-k-j}P)$
vanish for $0\le j<n-k$; thus the group 
$\HH{k}(Y,\Omega_Y^{n-k}(\log D\ichichi)\tensor L\inv)$ and the kernel of
the map 
$\HH{n}\left(Y,S^{n-k}(V^{\chi})^*\tensor
L\inv\right)\to\HH{n}\left(Y,S^{n-k-1}(V^{\chi})^*\tensor P^{\chi}\tensor
L\inv\right)$ are naturally isomorphic. By Serre duality, the latter group
is dual to
$R^{n-k,\chi}_{\omega_Y\tensor L}$.
\qed
The next result is the analogue in our setting of Macaulay's duality
theorem.
\begin{prop}\label{macaulay} 
Assume that the cover $\cover$ satisfies property $(A)$.
For
$\chi\in G^*\setminus\{1\}$, set\,
$\omega_{\chi}=\omega_Y^2(D\ichichi)$\,; then:

i) there is a natural isomorphism $R^{n,\chi}_{\omega_{\chi}}\cong \C$

ii) let $L$ be a line bundle on $Y$ such that $L$  and $L\inv(D\ichichi)$
 (or$L\tensor
\omega_Y\inv$ and $(\omega_Y\tensor
L)\inv(D\ichichi)$) belong to $\Gamma_f\setminus\{0\}$; then
the multiplication map
$R^{k,\chi}_{\omega_Y\tensor L}\otimes R^{n-k,\chi}_{\omega_Y\tensor
L\inv(D\ichichi)}\to R^{n,\chi}_{\omega_{\chi}}$ is a perfect pairing,
corresponding to Serre duality via the isomorphism of lemma \ref{iso}. In
particular, one has natural isomorphisms:
$$U^{k,\chi}\cong R^{k,\chi}_{\omega_Y\tensor L\ichinv}.$$
\end{prop}
\proof In order to prove i), consider the complex (\ref{reslog}) for
$k=n$: twisting it by $\omega_Y(D\ichichi)\inv$ and arguing as in the
proof of lemma \ref{iso}, one shows the existence of a natural
isomorphism between
$R^{n,\chi}_{\omega_{\chi}}$ and $\HH{0}(Y,\OY)=\C$. 
In order to prove statement ii), one remarks that  the
group
 $\HH{k}(Y,\Omega_Y^{n-k}(\log D\ichichi)\tensor L\inv)$ is Serre
dual to
$\HH{n-k}(Y,\Omega^k_Y(\log D\ichichi)\tensor L(-D\ichichi))$ by
(\ref{dualita'}). By lemma
\ref{iso}  the latter group  equals
$( R^{k,\chi}_{\omega_Y\tensor L\inv(D\ichichi)})^*$. Both these
isomorphisms and the multiplication map are natural, and therefore
compatible with Serre duality. The last claim follows in view of
(\ref{fundrelbis}).
\qed

\section{Infinitesimal Torelli}\label{infinitesimal}
\setcounter{defn}{0}
\setcounter{equation}{0}

In this section we exploit the algebraic description of the IVHS of
a $G$-cover to prove an infinitesimal Torelli theorem.  We
will use freely the notation introduced in section \ref{ivhs}.

We recall that  {\em infinitesimal Torelli}
for the periods of
$k$-forms holds for a variety  $X$ if the map (\ref{inftor}) is
injective. By the remarks at the beginning of section \ref{ivhs}, a
$G$-cover
$\cover$ satisfies  infinitesimal Torelli property if for each character
$\chi\in G^*$ the intersection, as
$\phi$ varies in $G^*$, of the kernels of the maps
$\rho^k\ichifi$ of (\ref{inftorab}) is equal to zero.
The next theorem shows that this is actually the case for $k=n$,
under some ampleness assumptions on the building data of $\cover$. 
\begin{thm}\label{mainthm1}
Let $X$, $Y$ be smooth complete algebraic varieties of dimension
$n\geq 2$ and let $\cover$ be a $G$-cover with building data $L\ichi$,
$D_i$, $\chi\in G^*\setminus\{1\}$, $i=1,\ldots r$. If properties $(A)$
and $(B)$ are satisfied, then the following map is injective:
$$\Ho^1(X,T_X)\to \Hom (\Ho^0(X,\omega_X),\Ho^1(X,\Omega_X^{n-1})),$$
and, as a consequence, infinitesimal Torelli for the periods of $n$-forms
holds for
$X$.
\end{thm}
Before giving the proof, we deduce two effective results from theorem
\ref{mainthm1}.
\begin{thm}\label{mainthm1p}
Let $f:X\to\PP^n$, $n\ge 2$, be a $G$-cover with building data $L\ichi$,
$D_i$, $\chi\in G^*\setminus\{1\}$, $i=1,\ldots r$. Assume that
$L\ichi\tensor\omega_{\PP^n}>0$ and
$D_i\tensor\omega_{\PP^n}>0$ for $\chi\in G^*\setminus\{1\}$,
$i=1,\ldots r$; then  infinitesimal Torelli for the periods of
$n$-forms holds for $X$.
\end{thm}
\proof: by remark \ref{projective}, properties
$(A)$ and $(B)$ are satisfied in this case. 
\qed
\begin{thm}\label{mainthm1e}
Let $X$, $Y$ be smooth complete algebraic varieties of dimension
$n\geq 2$ and let $\cover$ be a $G$-cover with building data $L\ichi$,
$D_i$, $\chi\in G^*\setminus\{1\}$, $i=1,\ldots r$. Let $E$ be a very
ample divisor on
$Y$ and let $E_n$ be defined as in prop. \ref{effective}; if 
$D_i, L\ichi, \omega_Y(D_i), \omega_Y\tensor
L\ichi>E_n$ and $L\ichi, D_i\ge(n+1)E$ for $\chi\in G^*\setminus\{1\}$,
$i=1,\ldots r$,\/then infinitesimal Torelli for the periods of $n$-forms
holds for
$X$.
\end{thm}
\proof
Follows from theorem \ref{mainthm1} together with proposition
\ref{effective}.
\qed

\noindent{\bf Proof of theorem \ref{mainthm1}:}

Since all cohomology groups appearing in this proof are computed on $Y$,
we will omit $Y$ from the notation.

\noindent By proposition 
\ref{projform} and by the discussion at the beginning of the section, we
have to show that for every 
$\chi\in G^*$ the intersection of the kernels of the maps $$\rho\ichifi:
\Ho^1(T_Y(-\log \Delta\ichi)\tensor
L\ichi\inv)\to\Hom(U^{0,\phi},U^{1,\chi\phi}),$$ as $\phi$ varies in
$G^*$, is equal to zero. For $\chi,\phi\in G^*$, set
$R\ichifi=D\ichifinv-D_{\chi\phi,\phi\inv}$. Notice that $R\ichifi$ is
effective.  By (\ref{dualita'}) and (\ref{fundrel}) there is a natural
identification:
$$\Omega_Y^{n-1}(\log D\ichifinv)\tensor(\omega_Y\tensor
L_{\chi\phi}\tensor L\ifinv)\inv = T_Y(-\log D\ichifinv)\tensor
L\ichi\inv (R\ichifi).$$  So
the map
$\rho\ichifi$ can be viewed as the composition of the map 
$$i\ichifi:\Ho^1(T_Y(-\log \Delta\ichi)\tensor L\ichi\inv) \to 
\Ho^1(T_Y(-\log D\ichifinv) \tensor L\ichi\inv (R\ichifi)),$$
induced by inclusion of sheaves, and of the map
$$r_{\chi\phi,\phi}:\Ho^1(\Omega_Y^{n-1}(\log
D\ichifinv)\tensor(\omega_Y\tensor L_{\chi\phi}\tensor L\ifinv)\inv)
\to \Hom(U^{0,\phi},U^{1,\chi\phi})$$
induced by cup-product. Arguing
as in the proof of thm. $3.1$ of \cite{torelli}, one can show that, for
fixed
$\chi\in G^*$, the intersection of the kernels of $i\ichifi$, as $\phi$
varies in $G^*\setminus\{1,\chi\inv\}$, is zero. (Notice that lemma $3.1$
of
\cite{torelli}, although stated for surfaces, actually holds for
varieties of any dimension, and that the ampleness assumptions on the
building data allow one to apply it, in view of thm. \ref{kodaira}.)
So the statement will follow if we prove that the map
$r\ichifi:\Ho^1(\Omega^{n-1}_Y(\log D\ichichi)\tensor(\omega_Y\tensor
L\ichi\tensor L\ifinv)\inv) \to$ $\Hom(U^{0,\phi},
U^{1,\chi})$ is injective
for every pair $\chi$, $\phi$ of nontrivial characters. By lemma
\ref{iso}, the map $r\ichifi$ may be rewritten as:
$r\ichifi:(R^{n-1,\chi}_{\omega_Y^2\tensor L\ichi\tensor L\ifinv})^*\to
(R^{0,\chi}_{\omega_Y\tensor L\ifinv})^*\otimes
(R^{n-1,\chi}_{\omega_Y\tensor L\ichi})^*$. We prove that
$r\ichifi$ is injective by showing that the dual map
$r\ichifi^*:R^{0,\chi}_{\omega_Y\tensor L\ifinv}\otimes
R^{n-1,\chi}_{\omega_Y\tensor L\ichi}\to R^{n-1,\chi}_{\omega_Y^2\tensor
L\ichi\tensor L\ifinv}$, induced by
multiplication, is surjective.   In order to do this, it is
sufficient to observe that the multiplication map
$$\HH{0}(Y,\omega_Y\tensor L\ifinv)\tensor
\HH{0}(Y,S^{n-1}(V^{\chi})\tensor\omega_Y\tensor
L\ichi)\to\HH{0}(Y,S^{n-1}(V^{\chi})\tensor\omega_Y^2\tensor
L\ifinv\tensor L\ichi)$$
is surjective by property $(B)$.
\qed

\begin{rem}
In section 6 of \cite{moduli}, it is proven that for any abelian group
$G$ there exist families of smooth $G$-covers with ample
canonical class that are {\em generically complete}. In those cases, thm.
\ref{mainthm1} means that the period map is \'etale on a whole component
of the moduli space.
\end{rem}

\section{Sufficiently ample line bundles}
\setcounter{defn}{0}
\setcounter{equation}{0}

We take up the following definition from \cite{suffampio}
\begin{defn}
A property is said to hold for a {\em sufficiently ample} line bundle $L$
on the smooth projective variety
$Y$ if there exists an ample line bundle $L_0$ such that the property
holds whenever the bundle $L\tensor L_0\inv$ is ample. We will denote
this by writing that the property holds for $L>>0$.
\end{defn}
In this section, we collect 
some facts about sufficiently ample line bundles that will be used
to prove our variational Torelli result. In particular, we prove a variant
of proposition
$5.1$ of
\cite{ciclico} to the effect that, given sufficiently ample line bundles
$L_1$ and $L_2$ on $Y$, it is possible to recover $Y$ from the kernel of
the multiplication map $\Ho^0(Y,L_1)\tensor
\Ho^0(Y,L_2)\to\Ho^0(Y,L_1\tensor L_2)$.

The next lemma is ``folklore''. The  proof given here has been
communicated to the author by Mark Green.
\begin{lem}\label{green} Let ${\cal F}$ be a coherent sheaf on $Y$ and
let $L$ be a line bundle on $Y$. Then, if $L>>0$,
$$\Ho^i(Y,\F\tensor L)=0, \quad i>0.$$
\end{lem}
\proof  We proceed by descending induction on  $i$. If $i>n$, then the
statement is evident. Otherwise, fix an ample divisor $E$ and an integer
$m$ such that $\F(mE)$ is generated by global sections. This gives rise
to an exact sequence: $0\to\F_1\to \oplus \OY(-mE)\to \F\to 0$. Tensoring
with
$L$ and considering the corresponding long cohomology sequence, one sees
that it is enough that $\Ho^i(Y,L(-mE))=\Ho^{i+1}(Y,\F_1\tensor L)=0$,
 for $L>>0$. The vanishing of the former group follows
from Kodaira vanishing and the vanishing of the latter follows from the
inductive hypothesis.
\qed

\begin{lem}\label{genproj} Let $|E|$ be a very ample lin\-ear
sys\-tem on the smooth
projective variety $Y$ of dimension $n$. Denote by ${\cal M}_y$ the ideal
sheaf of a point $y\in Y$: if $L$ is a sufficiently ample line bundle on
$Y$, then the multiplication map:
$$\Ho^0(Y,L)\tensor\Ho^0(Y,{\cal M}_y(E))\to\Ho^0(Y,{\cal M}_y(L+E))$$ is
surjective for every $y\in Y$.
\end{lem}
\proof For $y\in Y$, set $V_y=\Ho^0(Y,{\cal M}_y(E))$ and consider the
natural exact sequence $0\to N_y\to V_y\tensor\OY\to{\cal M}_y(E)\to 0$;
tensoring with $L$ and considering the corresponding cohomology sequence,
one sees that if $\Ho^1(Y,N_y\tensor L)=0$ then the map $\Ho^0(Y,L)\tensor
V_y\to\Ho^0(Y,{\cal M}_y(L+E))$ is surjective. By lemma \ref{green},
there exists an ample line bundle $L_y$ such that $\Ho^1(Y,N_y\tensor
L)=0$ if $L> L_y$.  In order to deduce from this the
existence of a line bundle $L_0$ such that $\Ho^1(Y,N_y\tensor L)=0$ for
every $y\in Y$
 if $L> L_0$, we proceed as
follows. Consider the product
$Y\times Y$, with projections
$p_i$,
$i=1,2$, denote by
$\I_{\Delta}$ the ideal sheaf of the diagonal in $Y\times Y$ and set
${\cal V}=p_{2*}(p_1^*\OY(E)\tensor \I_{\Delta})$. ${\cal V}$ is a fibre
bundle on
$Y$ such that the fibre of ${\cal V}$ at $y$ can be naturally identified
with $V_y$. We define the sheaf ${\cal N}$ on $Y\times Y$ to be the kernel
of the map
$p_2^*{\cal V}\to p_1^*\OY(E)\tensor \I_{\Delta}$; the restriction of
${\cal N}$ to
$p_2\inv(y)$ is precisely $N_y$. For any fixed line bundle $L$ on $Y$,
$h^1(Y,N_y\tensor L)=h^1(p_2\inv(y),{\cal N}\tensor p_1^*L|_{p_2\inv
(y)})$ is an upper-semicontinuous function of $y$. Thus, we may find a
finite open covering $U_1,\ldots U_k$ of $Y$ and ample line bundles
$L_1,\ldots L_k$ such that
$h^1(Y,N_y\tensor L)=0$ for $y\in U_i$ if $L\tensor L_i\inv$ is ample. To
finish the proof it is enough to set $L_0=L_1\tensor\ldots\tensor L_k$.
\qed

\begin{lem}\label{veryample} Let $Y$ be a smooth projective variety of
dimension $n\ge 2$, $Y\ne\PP^n$. If $E$ is a very ample divisor on
$Y$, then:

i) if
$L>\omega_Y((n-1)E)$, then $L$ is base point free.

 ii) if $L>\omega_Y(nE)$, then $L$ is very
ample.
\end{lem}
\proof In order to prove the claim, it is enough to show that if $C$ is a
smooth curve on $Y$ which is the intersection of $n-1$ divisors of
$|E|$, then $L|_C$ is  base point free (very ample) and $|L|$ restricts to
the complete linear $L|_C$. In view of the assumption, the former
statement follows from adjunction on
$Y$ and Riemann-Roch on $C$, and the latter follows from the vanishing of
$\Ho^1(Y,\I_C\tensor L)$, where $\I_C$ is the ideal sheaf of $C$. In
turn, this vanishing can be shown by means of the Koszul complex
resolution:
$$0\to\L((1-n)E)\to L\tensor\wedge^{n-2}\left(\oplus_1^{n-1}\OY
(-E)\right)\ldots\to
\oplus_1^{n-1}L(-E)\to \I_C\tensor L\to 0.$$  
\qed

The following proposition is a variant for sufficiently ample line
bundles of prop.
$5.1$ of
\cite{ciclico}. Our statement is weaker, but we do not need to assume
that one the line bundles involved  is much more ample than the other
one.
\begin{prop}\label{11}
  Let $L_1$ and $L_2$ be very ample line bundles on a smooth
projective variety $Y$. Let $\phi_1:Y\to\PP_1$ and $\phi_2:Y\to\PP_2$ be
the corresponding embeddings into projective space, and let
$f:Y\to\PP_1\times\PP_2$ be the composition of the diagonal embedding of
$Y$ in $Y\times Y$ with the product map $\phi_1\times\phi_2$. If
$L_1,L_2>>0$, then $f(Y)$ is the zero set of the elements of
$\Ho^0(\PP_1\times\PP_2,{\cal O}_{\PP_1\times\PP_2}(1,1))$ vanishing on
it.
\end{prop}
\proof By \cite{ciclico}, prop. 5.1, we may find very ample line bundles
$M_i$,
$i=1,2$ such that, if  $\psi_i:Y\to \Q_i$ are the corresponding
embeddings in projective space and $g:Y\to\Q_1\times\Q_2$ is the
composition of the diagonal embedding of $Y$ in $Y\times Y$ with
$\psi_1\times\psi_2$, then $g(Y)$ is the scheme-theoretic intersection of
the elements of $\Ho^0(\Q_1\times\Q_2,{\cal O}_{\Q_1\times\Q_2}(1,1))$
vanishing on it. By prop. \ref{veryample}, we may assume that $L_i\tensor
M_i\inv$ is very ample, $i=1,2$. To a divisor $D_i$ in $|L_i\tensor
M_i\inv|$ there corresponds a projection $p_{D_i}:\PP_i\cdots\to \Q_i$
such that
$\psi_i=p_{D_i}\circ \phi_i$. The claim will follow if we show that for
$(x_1,x_2)\notin f(Y)$ one can find $D_i\in |L_i\tensor M_i\inv|$ such
that $p_{D_i}$ is defined at $x_i$, $i=1,2$, and 
$(p_{D_1}(x_1), p_{D_2}(x_2))\notin g(Y)$. In fact, this implies that
there exists  $s\in \Ho^0(\Q_1\times \Q_2,{\cal
O}_{\Q_1\times\Q_2}(1,1))$ that vanishes on
$g(Y)$ and  does not vanish at $(p_{D_1}(x_1), p_{D_2}(x_2))$, so that the
pull-back of 
$s$ via $p_{D_1}\times p_{D_2}$ is a section of $\Ho^0(\PP_1\times
\PP_2,{\cal O}_{\PP_1\times\PP_2}(1,1))$ that vanishes on $f(Y)$ and does
not vanish at $(x_1,x_2)$.

By lemma
\ref{genproj}, if $L_i>>0$, then the multiplication
map
$\Ho^0(Y,L_i\tensor M_i\inv)\tensor\Ho^0(Y,M_i\tensor{\cal
M}_y)\to\Ho^0(Y,L_i\tensor{\cal M}_y)$ is surjective $\forall
y\in Y$ and for $i=1,2$. Notice that, in particular, this implies that
the map
$\Ho^0(Y,L_i\tensor M_i\inv)\tensor \Ho^0(Y,M_i)\to\HH{0}(Y,M_i\tensor
L_i)$ is surjective for
$i=1,2$, so that, for a generic choice of
$D_i$, the projection $p_{D_i}$ is defined at $x_i$. If either  $(x_1,x_2)\in
\phi_1(Y)\times\phi_2(Y)$ or there exists a
divisor $D_i$, for
$i=1$ or $i=2$, such that $p_{D_i}(x_i)\notin\psi_i(Y)$, then we are set. So
assume that, say, $x_1\notin \phi_1(Y)$ and that
$p_{D_i}(x_i)\in\psi_i(Y)$, for a generic choice of
$D_i$, for $i=1,2$. Fix a divisor $D_2$ such that
$p_{D_2}$ is defined at $x_2$ and write $p_{D_2}(x_2)=\psi_2(y_2)$, with
$y_2\in Y$. Now it is enough to show that there exists $D_1$ such that
$p_{D_1}$ is defined at $x_1$ and $p_{D_1}(x_1)\ne \psi_1(y_2)$. Since the
multiplication map
$\Ho^0(Y,L_1\tensor M_1\inv)\tensor\Ho^0(Y,M_1\tensor{\cal
M}_{y_2})\to\Ho^0(Y,L_1\tensor{\cal M}_{y_2})$ is surjective, there exist
$\sigma\in \Ho^0(Y,L_1\tensor M_1\inv)$ and $\tau \in
\Ho^0(Y,M_1\tensor{\cal M}_{y_2})$ such that $\sigma\tau$ corresponds to
a hyperplane of $\PP_1$ passing through $\phi_1(y_2)$ but not through
$x_1$. If
$D_1$ is the divisor of
$\sigma$, then the projection $p_{D_1}$ is defined at $x_1$ and
$p_{D_1}(x_1)\ne\psi_1(y_2)$.
\qed

\section{A variational Torelli theorem}
\setcounter{defn}{0}
\setcounter{equation}{0}

In this section we prove a variational Torelli theorem for the family
${\cal X}\to W$ of abelian covers with fixed basis and fixed $L\ichi$'s.

Let ${\cal Z}\to B$ be a flat family of smooth projective
polarized varieties on which $G$ acts fibrewise, and let $Z$ be the fibre
of ${\cal Z}$ over the point $0\in B$. It is possible to show that the
monodromy action on the cohomology of
$X$ preserves the group action;  therefore one may define a
{\em
$G$-period map}, by dividing the period domain $D$ by the subgroup of
linear transformations that, beside preserving the integral lattice and
the polarization, are compatible with the $G$-action. In particular, let
$\cover$ be a $G$-cover, with $X$ and $Y$  smooth projective varieties
of dimension
$n\ge 2$ and with building data $L\ichi$, $D_i$. In section \ref{covers},
we have introduced the family
${\cal X}\to W$ of
$G$-covers of $Y$, obtained by letting the sections
$s_i\in\HH{0}(Y,\OY(D_i))$ vary in equations (\ref{equazioni}). 
There is an obvious $G$-action on ${\cal X}$, and
the choice of an ample divisor $E$ on $Y$ gives a $G$-invariant
polarization of ${\cal X}$.  Since two
elements
$(s_1,\ldots s_r)$ and $(s'_1,\ldots s'_r)$ of $W$ represent the same
$G$-cover iff there exist $\lambda_i\in\C^*$ such that
$s'_i=\lambda_is_i$, the $G$-period map can be regarded as being defined
on the image $\tilde{W}$ of $W$ in
$|D_1|\times\cdots\times|D_r|$.
We denote by $s$ the image in $\tilde{W}$ of the point $(s_1,\ldots
s_r)$, by $X_s$ the corresponding $G$-cover of $Y$, by $T_s$ the
tangent space to $\tilde{W}$ at $s$ and, consistently with the notation
of section \ref{ivhs}, by $U_s^{k,\chi}$ the subspace of
$\HH{k}(X_s,\Omega^{n-k}_{X_s})$ on which $G$ acts via the character
$\chi$.  Remark that the space
$T_s$  can be naturally
identified with
$\oplus_i\HH{0}(Y,\OY(D_i))/(s_i)$. Denote by $\Gamma$ the subgroup of
$GL(T_s)\times
GL(\HH{0}(X_s,\omega_{X_s}))\times GL(\HH{1}(X_s,\Omega^{n-1}_{X_s}))$
that preserves the $G$-action on the cohomology of $X_s$.
\begin{thm}\label{mainthm2}
Assume that the dimension $n$ of $Y$ is $\ge2$, that the canonical class
$\omega_Y$ of
$Y$ is very ample and that $L\ichi$,
$D_i>>0$, $\chi\ne 1$, $i=1,\ldots r$. Assume that for $i=1,\ldots r$
the identity is the only automorphism of $Y$ that preserves the linear
equivalence class of $D_i$; moreover assume that 
$\forall i=1,\ldots r$ there exists $\chi\in G^*$ (possibly depending on $i$)
such that
$\chi(g_i)\ne 1$ and
$L\ichi>D_i$. Then a generic point
$s\in\tilde{W}$ is determined by the
$\Gamma$-class of the linear map:
$$T_s\to\oplus_{\chi}\Hom(U^{0,\chi}_s,U^{1,\chi}_s)),$$
which represents the first piece of the algebraic part of the IVHS for
$n$-forms.
\end{thm}
\begin{rem} The assumption, made in thm. \ref{mainthm2} and corollary
\ref{mainthm2bis}, that for every $i=1,\ldots r$ there exists $\chi\in
G^*$ such that $\chi(g_i)\ne 1$ and $L\ichi>D_i$ is not satisfied by
simple cyclic covers, namely totally ramified covers branched on an
irreducible divisor. Still there are many cases in which our results
apply: for instance, construction $6.2$ of \cite{moduli} provides
examples with
$G=\Z/_{m_1}\times\cdots\times\Z/_{m_{r-1}}$, $r-1\ge n$, branched on $r$
algebraically equivalent divisors $D_1,\dots D_r$. The ramification order
over $D_i$ is equal to
$m_i$ for
$i<r$, and it is equal to the least common multiple of the $m_i$'s for
$i=r$. (The assumption, made in \cite{moduli}, that $m_i|m_{i+1}$ is
actually unnecessary in order to make the construction.) Apart from the
case $r=m_1=m_2=2$, if the branch divisors are  ample, then
for every $i$ there exists $\chi$ such that $L\ichi(-D_i)$ is ample.
\end{rem}
Before giving the proof of
thm. \ref{mainthm2}, we state:
\begin{cor}\label{mainthm2bis}
Under the same assumptions as in theorem \ref{mainthm2} 
the $G$-period map for
$n$-forms has degree
$1$ on
$\tilde{W}$.
\end{cor}
\proof
Let ${\cal X}\to B$  be a family of
polarized varieties and let $X_b$ be the fibre of ${\cal X}$ over a
point $b\in B$ and assume that, if there is an isomorphism of the IVHS's
of
$X_b$ and
$X_{b'}$  preserving the polarization and the real structure, then the
varieties $X_b$ and $X_{b'}$ are isomorphic ( this property is usually
expressed by saying that ``variational Torelli holds''). In \cite{cdt} it
is proven that in this case the period map is generically injective on
$B$ up to isomorphism of varieties. Using exactly the same arguments, one
can show that if $G$ acts on the family ${\cal X}\to B$ fibrewise and
if the IVHS of a fibre $X_b$ determines $X_b$ up to isomorphisms
preserving the $G$-action, then the $G$-period map is generically
injective on $B$.
\qed
\noindent{\bf Proof of theorem \ref{mainthm2}:}
Whenever confusion is not likely to arise, we omit to write the space
where cohomology groups are computed.
 
By theorem $0.3$ of \cite{suffampio}, if
$D_i>>0$ for
$i=1,\ldots r$, then the period map $P$ for $n-1$ forms has degree $1$ on
$|D_i|$. Denote by $U_i$ the open subset of
$|D_i|$ consisting of the points $z$ such
that $P\inv(P(z))=\{z\}$, and fix $s\in
\tilde{W}\cap\left(U_1\times\ldots\times U_r\right)$. From now on we will
drop the subscript $s$ and write $X$ for $X_s$, $T$ for $T_s$, and so on.
For each index $i\in\{1,\ldots r\}$  set $S_i=\{\chi\in
G^*|\chi(g_i)=1\}$; as a first step, we show  that the subspace
$\HH{0}(\OY(D_i))/(s_i)$ of $T$ is the
intersection of the kernels of the maps $\rho\ichi:T\to
\Hom(U^{0,\chi},U^{1,\chi})$, as $\chi$ varies in $S_i\setminus\{1\}$.
Since
$\omega_Y$ is ample,  $\HH{0}(T_Y)=0$ by thm. \ref{kodaira}, and so $T$ equals
$R^{1,1}_{\OY}$. The map $\rho\ichi$ factors through the surjection
$R^{1,1}_{\OY}\to R^{1,\chi}_{\OY}$, whose kernel is
$\oplus_{\{i|\chi(g_i)= 1\}}\HH{0}(Y,\OY(D_i))/(s_i)$. In turn, by
sequence \ref{tlog}, $R^{1,\chi}_{\OY}$ injects in $\HH{1}(Y,T_Y(-\log
D\ichichi))$. We have shown in the proof of thm. \ref{mainthm1} that the
map $\HH{1}(T_Y(-\log D\ichichi))\to \Hom(U^{0,\chi},U^{1,\chi})$ is
injective. So $\ker\rho\ichi=\oplus_{\{i|\chi(g_i)=
1\}}\HH{0}(Y,\OY(D_i))/(s_i)$. As we have remarked in section
\ref{covers}, if $i\ne j$, then the subgroups of $G$ generated by $g_i$
and
$g_j$ intersect only in  $\{0\}$, and so $i$ is the only index such
that 
$\chi(g_i)=1$ for all $\chi\in S_i$. We conclude that
$\cap_{\chi\in S_i\setminus\{1\}}\ker\rho\ichi=\HH{0}(Y,\OY(D_i))/(s_i)$.

Now fix $i\in\{1,\ldots r\}$ and let $\chi\in G^*\setminus S_i$  be such
that $L\ichi>D_i$: the restriction to $\HH{0}(\OY(D_i))/(s_i)$ of the
map $T\tensor U^{0,\chi\inv}\to
U^{1,\chi\inv}$ is the multiplication map
$\HH{0}(\OY(D_i))/(s_i)\otimes\HH{0}(\omega_Y\tensor
L\ichinv)\to\HH{0}(Y,\omega_Y\tensor L\ichinv(D_i))|_{D_i}$, followed by
the inclusion $\HH{0}(\omega_Y\tensor
L\ichinv(D_i))|_{D_i}\to U^{1,\chi\inv}$.

{\em Claim:} the kernel of the latter map is
equal to zero.

If we assume that the claim holds, then we have
recovered the multiplication map
$\HH{0}(\OY(D_i))/(s_i)\otimes \HH{0}(\omega_Y\tensor
L\ichinv)\to \HH{0}(Y,\omega_Y\tensor
L\ichinv(D_i))|_{D_i}$. The right
kernel of this map is $(s_i)\HH{0}(\omega_Y\tensor L\ichinv)(-D_i))$. So
we can reconstruct the map $\HH{0}(\OY(D_i))|_{D_i}\otimes
\HH{0}(\omega_Y\tensor L\ichinv)|_{D_i}\to \HH{0}(\omega_Y\tensor
L\ichinv(D_i))|_{D_i}$. Let
$\phi_{1}:Y\to\PP_1$ be the embedding defined by the linear system
$|D_i|$, let
$\phi_{2}:Y\to\PP_2$ be the embedding defined by the linear system
$|\omega_Y\tensor L\ichinv|$  and let
$f:Y\to\PP_1\times\PP_2$ be the composition of the diagonal embedding
$Y\to Y\times Y$ with the product map $\phi_1\times\phi_2$: by prop.
\ref{11}, $f(Y)$ is the zero set of the elements of the
kernel of
$\HH{0}(\OY(D_i))\otimes\HH{0}(\omega_Y\tensor L\ichinv)\to
\HH{0}(\omega_Y\tensor L\ichinv(D_i))$. 
This implies that $f(D_i)$ the zero set in
$\PP\left(\HH{0}(\OY(D_i))|_{D_i}\right)\times
\PP\left(\HH{0}(\omega_Y\tensor L\ichinv)|_{D_i}\right)$ of the elements
of the kernel of 
 $\HH{0}(\OY(D_i))|_{D_i}\otimes
\HH{0}(\omega_Y\tensor L\ichinv)|_{D_i}\to \HH{0}(\omega_Y\tensor
L\ichinv(D_i))|_{D_i}$. Thus it is
possible to recover  $D_i$ as an abstract variety, for every $i=1,\ldots
r$. Since $s\in U_1\times\cdots\times U_r$, this is enough to determine
the
$D_i$'s as  divisors on $Y$ and, in turn, the point $s$.

In order to complete the proof, we have to prove the claim.  Let
$D_{i_1},\ldots D_{i_s}$ be the components of $D\ichichi$. (Recall that
there exists $j_0$ such that $i=i_{j_0}$.)
Let $V=\oplus_j\OY(D_{i_j})$, let $V_i=\oplus_{j\ne
j_0}\OY(D_{i_j})$, let
$P$ be the generalized prolongation bundle associated to $(D_{i_1},
\ldots D_{i_s})$ and let $P_i$ be the generalized prolongation
bundle associated to
$(D_{i_1},\ldots\hat{D_i},\ldots D_{i_s})$. There is a natural short exact
sequence $0\to P_i\to P\to\OY\to 0$, with dual sequence
$0\to\OY\to P^*\to P_i^*\to 0$. From this and sequence
\ref{tlog}, tensoring with
$\omega_Y\tensor L\ichinv$ and taking global sections, one deduces the
following commutative diagram with exact rows:
$$\begin{array}{ccccccccc}
0 &\!\!\!\!\!\! \to\!\!\!\!\!\! & \HH{0}(\omega_Y\tensor
L_{\chi\!\inv}\!) &
\!\!\!\!\!\!
\to\!\!\!\!\!\! &
\HH{0}\!(P^*\tensor\omega_Y\tensor L_{\chi\!\inv}\!)
&\!\!\!\!\!\!\to\!\!\!\!\!\! &
\HH{0}\!(P_i^*\tensor\omega_Y\tensor L_{\chi\!\inv}\!) &\!\!\!\!\!\!
\to\!\!\!\!\!\! & 0\\
\phantom{1} & \phantom{1} &\downarrow &\phantom{1} & \downarrow &
\phantom{1}& \downarrow &
\phantom{1}& \phantom{1}\\
0 & \!\!\!\!\!\! \to\!\!\!\!\!\! & \HH{0}\!(\omega_Y\tensor
L_{\chi\!\inv}\!(D_i)) &
\!\!\!\!\!\! \to\!\!\!\!\!\! &
\HH{0}\!(V\tensor\omega_Y\tensor L_{\chi\!\inv}\!)
&\!\!\!\!\!\!\to\!\!\!\!\!\! &
\HH{0}\!(V_i\tensor\omega_Y\tensor L_{\chi\!\inv}\!) & \!\!\!\!\!\!
\to\!\!\!\!\!\! & 0
\end{array}$$
In view of sequence \ref{tlog}, by applying snake's lemma to this diagram
one obtains the following exact sequence:
$\HH{0}(T_Y(-\log(D\ichichi-D_i))\tensor \omega_Y\tensor L\ichinv)\to
\HH{0}(\omega_Y\tensor L\ichinv(D_i))|_{D_i}\to U^{1,\chi\inv}$. So it is
enough to show that
$\HH{0}(T_Y(-\log(D\ichichi-D_i)\tensor L\ichinv)=0$. Using the
isomorphism (\ref{dualita'}) and the relations (\ref{fundrelbis}), one has
$T_Y(-\log(D\ichichi-D_i))\tensor\omega_Y\tensor L\ichinv\cong
\Omega^{n-1}_Y(\log(D\ichichi-D_i)\tensor L\ichi\inv(D_i)$.
In view of the assumptions, the required vanishing now follows from thm.
\ref{kodaira}.
\qed

\bigskip

Dipartimento di Matematica e Informatica -- Universit\`a di Udine

Via delle Scienze 206 - 33100 Udine (Italia)

E-mail: pardini@ten.dimi.uniud.it


\begin{thebibliography}{10}
\bibitem[1]{cdt} D.~Cox, R.~Donagi, L.~Tu, {\em Variational Torelli
implies generic Torelli}, Invent. Math. {\bf 88} (1987), 439-446.
\bibitem[2]{cggh} J.~Carlson, P.~Griffiths, {\em Infinitesimal
variations of Hodge structure (I)}, Compositio Mathematica, {\bf 50}
(1983), 109-205.
\bibitem[3]{donagi} R.~Donagi,{\em Generic Torelli for projective
hypersurfaces}, Compositio Mathematica, {\bf 50} (1983),
325-353.
\bibitem[4]{EV} H.~Esnault, E.~Viehweg, {\em Lectures on
Vanishing Theorems}, DMV Seminar, Band 20, Birkh\"auser 1992.
\bibitem[5]{el} L.~Ein, R.~Lazarsfeld {\em Syzygies and Koszul
cohomology of smooth projective varieties of arbitrary dimension},
Invent. Math., {\bf 111}, 51-67 (1993).
\bibitem[6]{moduli} B.~Fantechi, R.~Pardini, {\em Automorphisms and
moduli spaces of varieties with ample canonical class via deformations of
abelian covers}, preprint Max Planck Institut f\"ur Mathematik n. 112,
1994. 
\bibitem[7]{cime} M.~Green, J.~Murre, C.~Voisin, {\em Algebraic Cycles
and Hodge Theory}, C.I.M.E. notes Torino 1993, Springer L.M.N. {\bf
1594}.
\bibitem[8]{suffampio} M.~Green, {\em The period map for hypersurface
sections of high degree of an arbitrary variety}, Compositio Mathematica,
{\bf 55} (1984), 135--156.
\bibitem[9]{Har} R.~Hartshorne, {\em Algebraic Geometry},  GTM 52,
Springer 1977.
\bibitem[10]{ciclico} K.~Ivinskis, {\em A variational Torelli theorem
for cyclic coverings of high degree}, Compositio Mathematica, {\bf 85}
(1993), 201-228.
\bibitem[11]{konnociclico} K.~Konno, {\em On deformations and the local
Torelli problem of cyclic coverings}, Math. Ann. {\bf 231} (1977),
601-617.
\bibitem[12]{konnovar} K.~Konno, {\em On the variational Torelli problem
for complete intersections}, Compositio Mathematica, {\bf 78} (1991),
271-296.
\bibitem[13]{konnocompl} K.~Konno {\em Infinitesimal Torelli for
complete intersections in certain homogeneous K\"ahler manifolds},
Tohoku Math. J.,{\bf 38} (1986), 609-624.
\bibitem[14]{lwp} D.~Lieberman, R.~Wilsker, C.~Peters {\em A theorem of
local Torelli type}, Math. Ann., {\bf 231} (1977), 39-45.
\bibitem[15]{abeliani} R.~Pardini, {\em Abelian covers of algebraic
varieties}, J.\ reine angew.\ Math., {\bf 417} (1991), 191--213.
\bibitem[16]{torelli} R.~Pardini, {\em Infinitesimal Torelli and abelian
covers of algebraic surfaces}, in ``Problems in the theory of surfaces
and their classification", F.~Catanese, C.~Ciliberto and M.~Cornalba
eds., Symp.~Math.
\hbox{INDAM} XXXII, Academic Press, 1991.
\bibitem[17]{peters} C.~Peters {\em The local Torelli theorem II. Cyclic
branched coverings}, Ann. Scuola Normale Sup. Pisa Cl. Sci., (4) {\bf 3}
(1976), 321-340.
\bibitem[18]{SS} B.~Schiffman, A.~Sommese, {\em Vanishing Theorems on
Complex Manifolds}, Progr. in Math. {\bf 56}, Birkh\"auser (1985).
\end{thebibliography}
\end{document}